\documentclass[12pt]{article}
\usepackage{a4wide}
\usepackage{theorem}
\usepackage{amssymb}
\usepackage{amsmath}
\usepackage{amsfonts}
                     
\theoremstyle{plain}
\newtheorem{theo}{Theorem}
\newtheorem{proper}{Property}
\newtheorem{cons}{Consequence}[section]
\newtheorem{nota}[cons]{Notation}
\newtheorem{lem}[cons]{Lemma}
\newtheorem{rem}[cons]{Remark}
\newtheorem{defi}[cons]{Definition}
\newtheorem{coro}[cons]{Corollary}

\newcommand{\biindice}[3]%
{\renewcommand{\arraystretch}{0.5}
\begin{array}[t]{c}
#1\\
{\scriptscriptstyle #2}\\
{\scriptscriptstyle #3}
\end{array}
\renewcommand{\arraystretch}{1}}

\newcommand{\geqs}[1]{\arraycolsep0.1pt \renewcommand{\arraystretch}{0.5}
\begin{array}[t]{c}
\geq\\
#1
\end{array}
\arraycolsep5pt \renewcommand{\arraystretch}{1}
}

\newcommand{\leqs}[1]{\arraycolsep0.1pt \renewcommand{\arraystretch}{0.5}
\begin{array}[t]{c}
\leq\\
#1
\end{array}
\arraycolsep5pt \renewcommand{\arraystretch}{1}
}

\newcommand{\dN}{I \! \! N}
\newcommand{\dR}{I \! \! R}

\title{Fu\v cik Spectrum for the Neumann Problem with Indefinite Weights}
\author{~~\\[10mm]Mohssine ALIF \footnote{e-mail: malif@ictp.trieste.it}\\
D\'ep. Math., C.P. 214, Universit\'e Libre de Bruxelles,\\ 
1050 Bruxelles, Belgium\\
and\\
The Abdus Salam International Centre for Theoretical Physics,\\
Trieste, Italy
}
\date{}

\begin{document}

\maketitle

%\selectlanguage{french}

\begin{abstract}
We obtain a description of the Fu\v cik spectrum associated to the
one-dimensional asymmetric problem with indefinite weights
$Lu = am(t)u^{+} - bn(t)u^{-}$ in $]T_{1},T_{2}[$,
$u'(T_{1}) = 0 = u'(T_{2})$, where $L$ is a Sturm-Liouville operator.
Our approach is based on the shooting method.
\end{abstract}

%\selectlanguage{french}

\section{Introduction}

This paper is concerned with the study of the Fu\v cik spectrum $\Sigma$ 
associated to the semilinear Neumann problem with weights 
\begin{equation}
    \left\{
    \begin{array}{l}
        Lu = am(t)u^{+} - bn(t)u^{-} \mbox{ in } ]T_{1},T_{2}[,\\
        u'(T_{1}) = 0 = u'(T_{2}) ,
    \end{array}
\right.
    \label{eq:1.1}
\end{equation} 
where $Lu:=-[p(t)u']'+q(t)u$, $p,q,m$ and $n \in C[T_{1},T_{2}]$, $p(t) > 
0$ on $[T_{1},T_{2}]$, $q(t) \geq 0$ on $[T_{1},T_{2}]$, $m(t)$ and
$n(t)$ are both $\not \equiv 0$ and $u^{\pm}:=\max\{\pm
u,0\}$. $\Sigma$ is defined as the set of those  $(a,b) \in \dR^{2}$
such that the problem (\ref{eq:1.1}) has a nontrivial solution.
This set plays an important role in the
study of semilinear problems of type
\begin{equation}
    \left\{
    \begin{array}{l}
        Lu = f(t , u(t)) \mbox{ in } ]T_{1},T_{2}[,\\
	u'(T_{1}) = 0 = u'(T_{2}).
    \end{array}
\right.
    \label{eq:1.2}
\end{equation}
\par
Consider first the case where there is only one weight-function with a
constant sign in the whole of the interval $[T_1 , T_2]$, i.e. the case where
$m(t) = n(t)$ and for instance $m(t) > 0$ on $[T_1,T_2]$.
The problem (\ref{eq:1.1}) becomes:
\begin{equation}
    \left\{
    \begin{array}{l}
        Lu = m(t)[au^{+} - bu^{-}] \mbox{ in } ]T_{1},T_{2}[,\\
	u'(T_{1}) = 0 = u'(T_{2}).
    \end{array}
\right.
    \label{eq:1.3}
\end{equation}
It is well known (see \cite{Ca}, \cite{Dr}, \cite{Ry})
that, in this case, $\Sigma$ is made of the two lines
$\dR \times \{\lambda_{1}^{m}\}$ and $\{\lambda_{1}^{m}\} \times
\dR$ together with a sequence of hyperbolic like curves in $\dR^{+}\times
\dR^{+}$ passing through $(\lambda_{k}^{m}, \lambda_{k}^{m}$), $k \geq 2$;
each point $(\lambda_{k}^{m}, \lambda_{k}^{m})$ belongs to one (or
two) of these curves. Along this (or these two) curve(s), the
corresponding solutions for (\ref{eq:1.3})
have exactly $(k-1)$-zeros in $]T_{1},T_{2}[$. Here
$(0 =) \lambda_{1}^{m} < \lambda_{2}^{m} \leq \ldots \to + \infty$
designate the sequence of eigenvalues of the linear problem
associated to (\ref{eq:1.3}) :
\begin{equation}
    \left\{
    \begin{array}{l}
    Lu  = \lambda m(t) u \mbox{ in } ]T_{1},T_{2}[,\\
    u'(T_{1}) = 0 = u'(T_{2}).
\end{array}
\right.
    \label{eq:1.4}
\end{equation}
The description of  $\Sigma$ is explicit in the classical case where
$Lu = -u''$ and $m(t) \equiv n(t) = 1$
(see \cite{Da}, \cite{Dr}, \cite{Fu}). In the second paragraph of section 3,
we investigate the situation where the weight-function
$m(t)$ changes sign in the interval $[T_1,T_2]$, i.e. when both $m^+$ and
$m^-$ are $\not\equiv 0$ in $[T_1,T_2]$. In
this case, we prove that the description of $\Sigma$ is rather
comparable with that of the Fu\v cik spectrum associated to the Dirichlet
problem
\begin{equation}
    \left\{
    \begin{array}{l}
        Lu = m(t)[au^{+} - bu^{-}] \mbox{ in } ]T_{1},T_{2}[,\\
        u(T_{1}) = 0 = u(T_{2})
    \end{array}
\right.
    \label{eq:1.5}
\end{equation}
which was studied recently in \cite{Al1} and \cite{Al2}:
$\Sigma$ is made of the four lines
$\dR \times \{\lambda_{1}^{m}\}$,
$\{\lambda_{1}^{m}\} \times \dR$, $\dR \times \{\lambda_{-1}^{m}\}$ and
$\{\lambda_{-1}^{m}\} \times \dR$ together with a double sequence of
hyperbolic like curves in both 
$\dR^{+} \times \dR^{+}$ and $\dR^{-} \times \dR^{-}$ passing through
$(\lambda_k^m , \lambda_k^m)$ and $(\lambda_{-k}^m , \lambda_{-k}^m)$,
$k \geq 2$, and, moreover, a (non-zero) number of
additional curves which appear in the other two quadrants $\dR^{+}
\times \dR^{-}$ and $\dR^{-} \times \dR^{+}$. This number depends on
``the number of changes of sign'' of $m(t)$ in the interval $[T_1 ,
T_2]$. More precisely, if $m(t)$ ``changes sign $N (= 1, 2, ... , +\infty)$-times'' in
$]T_1 , T_2[$, then $\Sigma$ contains exactly
$(2N-1)$-hyperbolic like curves in $\dR^{+} \times \dR^{-}$ and also in
$\dR^{-} \times \dR^{+}$. These additional curves can be classified
according to the number of zeros of the corresponding solutions of the
problem (\ref{eq:1.3}). Here as below $- \infty \leftarrow \ldots \leq \lambda_{-2}^{m}
< \lambda_{-1}^{m}\; (\leq 0 \leq)\; \lambda_{1}^{m} < \lambda_{2}^{m} \leq
\ldots \to + \infty$ are the eigenvalues of the linear problem (\ref{eq:1.4}).
\vspace{3mm}
\par
Let us return to the problem (\ref{eq:1.1}) and suppose that
$m(t)$ and $n(t)$ both change sign in the interval $[T_1 ,
T_2]$. As it was done for the Dirichlet case in \cite{Al1} and \cite{Al2},
we prove in the first paragraph of section 3 that, beside the
trivial part consisting of the four lines
$\{\lambda_{1}^{m}\} \times \dR$, $\dR \times \{\lambda_{1}^{n}\}$,
$\{\lambda_{-1}^{m}\} \times \dR$ and $\dR \times \{\lambda_{-1}^{n}\}$,
the Fu\v cik spectrum $\Sigma$ associated to (\ref{eq:1.1}) is made in each
quadrant of $\dR^2$ of {\it a (non-zero) odd or infinite number} of
hyperbolic like curves. Moreover, $p, q, r, s$ being in $\dN$, we
prove the existence of two weight-functions $m(t)$ and $n(t)$ such that
the spectrum $\Sigma$ associated to (\ref{eq:1.1}) exactly contains 
$(2p+1), (2q+1), (2r+1)$ and $(2s+1)$
hyperbolic like curves in $\dR^{+} \times \dR^{+}$, $\dR^{-} \times
\dR^{-}$, $\dR^{+} \times \dR^{-}$ and $\dR^{-} \times \dR^{+}$
respectively. Note that, here as above, certain curves may be double and
are then counted for two.
\vspace{3mm}
\par
Section 4 is mainly devoted to the study of the 
asymptotic behaviour of the ``first curves'' of $\Sigma$ in
each quadrant of $\dR^2$, i.e. those which lie the closest to the
trivial horizontal and vertical lines. In the Dirichlet case (\cite{Al1},
\cite{Al2}), it was observed that for instance the first curve in
$\dR^+ \times \dR^+$ is non-asymptotic on any side to the trivial horizontal
and vertical lines {\it if and only if} both $m^+$ and $n^+$ have compact support
in $[T_1 , T_2]$. In our Neumann case, it turns out that this
non-asymptotic character holds always true for the first curve in
$\dR^+ \times \dR^+$ even if $m^+$ or $n^+$ do not have compact support in
$[T_1 , T_2]$.
\vspace{3mm}
\par
Section 2, which has a preliminary character, deals with the linear
equation
\begin{equation}
        Lu = am(t)u \mbox{ in } \dR  .
   \label{eq:1.6}
\end{equation}
We introduce, among others, two ``zero-functions''.
The first one sends any zero of a nontrivial solution $u$ of
(\ref{eq:1.6}) onto ``the following zero'' of $u$; the second one
sends any $a \in \dR$ onto ``the first zero'' $> T_1$ of a nontrivial
solution $u$ of (\ref{eq:1.6}) satisfying $u'(T_1) = 0$. We investigate
several properties of these functions.

\section{Zero-Functions}
\setcounter{equation}{0}

Consider the linear equation (\ref{eq:1.6}) where $L$ is defined as in the 
introduction and $m \in C[T_{1}, T_{2}]$, $m(t) \not \equiv 0$. 
From now on, in this section, it will be supposed that $p(t)$, $q(t)$ and
the weight-function $m(t)$ are extended from
$[T_1 , T_2]$ to the whole of $\dR$ preserving the
continuity and the fact that
$p_{1} \leq p(t) \leq p_{2}$ and $0 \leq q(t) \leq q_{2}$ for some constants
$p_{1} > 0$, $p_{2}$ and $q_{2}$.\\
From standard results on the ordinary differential equations (cf.
\cite{Co} or \cite{Ha}), it follows that, for each $s \in \dR$, there exists a 
unique solution $u(t)  = u(t;a,s)$, $u(s) = 0$ and $u'(s) = 1$, of the problem
(\ref{eq:1.6}).\\ 
Similarly, we have the existence and the unicity of $v(t) = v(t;a,s)$ : solution
of (\ref{eq:1.6}) such that $v'(s)=0$ and $v(s)= 1$.\\
Moreover, $u(t;a,s)$ and $v(t;a,s)$ are both $C^1$-functions of 
$(t,a,s)$.
\begin{defi}\label{defi:2.1}
    {\rm We define the three zero-functions
    $\varphi$, $\psi_1$ and $\psi_2$ as follows
    $$
    \begin{array}{rcl}
    \varphi_{a}(s) & : = & \min \{ t > s : u(t;a,s)=0\}\, ,\\
    \psi_{1}(a) & : = & \min \{ t > T_1 : v(t;a,T_1)=0\}\, ,\\
    \psi_{2}(a) & : = & \sup \{ t < T_2 : v(t;a,T_2)=0\}\, ,
    \end{array}
    $$
    for each $(a,s)\in \dR^2$,
    with $\varphi_a(s)$ (resp. $\psi_{1}(a)$) = $+\infty$ if
    $u(t;a,s)$ (resp. $v(t;a,T_1)$) does not have any zero $> s$ (resp. $T_1$). 
    Similarly, $\psi_{2}(a)$ = $-\infty$ if $v(t;a,T_2)$ does not vanish at any 
    point $t < T_2$.}
\end{defi}

Note that $\varphi_{a}(s)$ (resp. $\psi_{1}(a)$) is in fact the first zero 
following $s$ (resp. $T_1$) of, not only $u(.;a,s)$ (resp. $v(.;a,T_1)$),
but also of any nontrivial solution $u$ of (\ref{eq:1.6}) satisfying
the initial condition $u(s) = 0$ (resp. $u'(T_1) = 0$).\\
Similarly, $\psi_2(a)$ is the last zero $< T_2$
of any solution $u \not \equiv 0$ of (\ref{eq:1.6})
satisfying $u'(T_2) = 0$.\\
On the other hand, it follows from classical ODE results that the zeros of both 
$u$ and $v$ are isolated. Thus the above definitions are meaningful.
\begin{nota}
\label{nota:2.2} 
{\rm Let
\vspace{2mm}\\
\begin{tabular}{ccl}
$A$ & := & $\{(a,s) \in \dR^{2} : \varphi_{a}(s) < + \infty\}$,\\
$B_1$ & := & $\{a \in \dR : \psi_{1}(a) < + \infty\}$,\\
$B_2$ & := & $\{a \in \dR : \psi_{2}(a) > - \infty\}$
\end{tabular}
\vspace{2mm}\\
and, for each $s \in \dR$,
\vspace{2mm}\\
\begin{tabular}{ccl}
$\alpha_s^>$ & := & $\inf \{ t>s : m(t) > 0 \}$,\\ 
$\alpha_s^<$ & := & $\inf \{ t>s : m(t) < 0 \}$
\end{tabular}
\vspace{2mm}\\
with the convention: $\inf \emptyset = +\infty$.\\
In the simplest cases, $\alpha_{s}^{>}$ (resp. $\alpha_{s}^{<}$)
is the lower bound of the first positive (resp.
negative) bump of $m(t)$ situated at the right of $s$.
}
\end{nota}

The zero-function $\varphi$ has the following
\vspace{3mm}
\\
{\bf Properties} (see \cite{Al1} or \cite{Al2})\\
{\it
\begin{tabular}{rrl}
1) & (i) & $A$ is open,\\
 & (ii) & $\varphi$ : $A \longmapsto \dR$ is a $C^1$-function.
\vspace{2mm}\\
2) & (i) & $\forall a\in \dR$, $\varphi_{a}(s)$ is increasing with respect to
$s$, strictly in $A$,\\
 & (ii) & $\forall s\in \dR$, $\varphi_{a}(s)$ is decreasing with respect to
$a$ for $a \geq 0$, strictly in $A$,\\
 & (iii) & $\forall s\in \dR$, $\varphi_{a}(s)$ is increasing with respect to
$a$ for $a \leq 0$, strictly in $A$.
\vspace{2mm}\\
3) & (i) & $\partial \varphi_{a}(s)/\partial s  > 0$ for $(a,s) \in A$,\\ 
 & (ii) & $\partial \varphi_a(s)/\partial a < 0$ for $(a,s)  \in A$ 
and $a > 0$,\\
 & (iii) & $\partial \varphi_a(s)/\partial a > 0$ for
$(a,s) \in A$ and $a < 0$.
\vspace{2mm}\\
4) & (i) & $\forall s \in \dR$, $\lim\limits_{a \longrightarrow +\infty}
\varphi_{a}(s) = \alpha_{s}^{>}$,\\
 & (ii) & $\forall s \in \dR$, $\lim\limits_{a \longrightarrow -\infty}
\varphi_{a}(s) = \alpha_{s}^{<}$.
\vspace{2mm}\\
5) & & $\forall s \in \dR$, $\lim\limits_{a \longrightarrow 0}
\varphi_{a}(s) = +\infty$.
\end{tabular}
}
\begin{rem}
\label{rem:2.3}
{\rm We will restrict ourselves below, in this section, to the study of 
$\psi_1$. By the change of variable ``$\tilde t = -t$'', one can easily 
deduce the properties of $\psi_2$ from those of $\psi_1$.
}
\end{rem}

Indeed, let us denote $\tilde v(\tilde t) := v(-\tilde t) = v(t)$.
By this change of variable, our equation (\ref{eq:1.6}) becomes: 
$-[p(-\tilde t){\tilde v}']' + q(-\tilde t){\tilde v}
= am(-\tilde t)\tilde v$.\\
Denote by $\tilde \varphi$, $\tilde \psi_1$ and $\tilde \psi_2$,
respectively, the three zero-functions associated to the latter equation. 
It is clear that 
$\tilde v(\tilde t,a,\tilde s) = -v(-\tilde t,a,-\tilde s)$, 
$\tilde \psi_1(a) = -\psi_2(a)$ and 
$\tilde \psi_2(a) = -\psi_1(a).$\\
Hence, the conditions on $v(.,a,T_2)$ at $T_2$ become initial conditions
for $\tilde v(.,a,-T_2)$ at $-T_2$ and $\psi_2(a) = -\tilde \psi_1(a)$. Q. E. D.
\vspace{3mm}
\par 
Now, we'll see that the properties of $\psi_1$ are similar to those of
$\varphi$.
\begin{proper}
\label{proper:1}
(i) $B_1$ is an open subset of $\dR$,\\
(ii) $\psi_1 : B_1 \longrightarrow \dR$ is a $C^1$-function.
\end{proper}
 
{\bf Proof.}
Let $a_0 \in B_1$. Note first that
$\frac{\partial}{\partial t} v(t;a_0,T_1)/_{t=\psi_1(a_0)}$ can not be
$= 0$. This is a direct consequence of standard uniqueness theorems for
the ODE
(see, e.g. \cite{Co} or \cite{Ha}) and the fact that $v(\psi_1(a_0);a_0,T_1) = 0$ and
$v(.;a_0,T_1) \not\equiv 0$. Hence using the implicit function theorem,
one gets open neighbourhoods $U$ of $a_0$ and $V$ of $\psi_1(a_0)$ and
a $C^1$-function $\tilde \psi_1 : U \longrightarrow V$ such that
$$\tilde \psi_1(a) = t \Longleftrightarrow (t,a) \in V \times U \;
\mbox{and} \; v(t;a,T_1) = 0 \, .$$
Now by the same arguments used in \cite{Al2} to prove property 1 of
$\phi$,
we can show that $\tilde \psi_1 \equiv \psi_1$ near $a_0$ and that
$U \subset B_1$, $U$ being sufficiently reduced. Q. E. D.
\vspace{3mm}
\par
The following two properties are concerned with the monotonicity and the
regularity of $\psi_1$.
\begin{proper}
\label{proper:2}
(i) $\psi_1$ is $\searrow$ in $\dR^+$, strictly in
$B_1 \cap \dR^+$,\\
(ii) $\psi_1$ is $\nearrow$ in $\dR^-$, strictly in
$B_1 \cap \dR^-$.
\end{proper}

{\bf Proof.} 
We'll prove only the first assertion (the proof of the second one is
similar). So let $a \in B_1 \cap
\dR^+$ and $\tilde a > a$ and let us verify that
$\psi_1(\tilde a) < \psi_1(a)$.
This will be clearly done if we succeed in proving that the solution
${\tilde v}(.,{\tilde a},T_1)$ of the linear equation 
\begin{equation}
L\tilde v = \tilde am(t)\tilde v 
\label{eq:2.1}
\end{equation}
has at least one zero in the interval $]T_1 , \psi_1(a)[$.\\
Writing the equations (\ref{eq:1.6}) and (\ref{eq:2.1}), respectively, as
$$
-\frac {1}{a}[p(t)v']' + \frac{1}{a}q(t)v = m(t)v
$$
and
$$
-\frac {1}{\tilde a}[p(t){\tilde v}']' +
\frac {1}{\tilde a}q(t){\tilde v} = m(t)\tilde v
$$
and applying the following lemma to these two latters on the interval
$[T_1 , \psi_1(a)]$, one gets the desired result.\\
By the same arguments, we prove that if $\tilde a \geq a \geq 0$, then
$\psi_1(\tilde a) \leq \psi_1(a)$. Property \ref{proper:2} is proved. Q. E. D.
\vspace{3mm}
\\
The following lemma will be used below repeatedly.
\begin{lem}
\label{lem:2.4}
Consider the following two equations:
$$
-[p_i(t){v_i}']' + q_i(t)v_i = m_i(t)v_i , \,\, t \in \dR,
\leqno{(*)_i}
$$
($i = 1,2$) and let $v_1$ and $v_2$ be two solutions of $(*)_1$ and $(*)_2$,
respectively,
such that
$$ \left\{ \begin{array}{rcl} v_i(T_1) & = & 1\,\,\, ,\\ 
{v_i}'(T_1) & = & 0\,\, (i = 1,2)\, .\end{array} \right. $$
Suppose that $p_1 \geq p_2\, (> 0)$, $q_1 \geq q_2\, (\geq 0)$, $m_1 \leq
m_2$ and that $v_1$ vanishes at a point $t_0 > T_1$. Then $v_2$ has at
least one zero in the interval $]T_1 , t_0]$.\\
If, moreover, at least one of the above three inequalities is strict, then
$v_2$ has a zero in the open interval $]T_1 , t_0[$.
\end{lem}

{\bf Proof.} 
Suppose $p_1 \geq p_2\, > 0$, $q_1 \geq q_2\, \geq 0$, $m_1 \leq m_2$
and suppose, by contradiction, that $v_1(t_0) = 0$ and $v_2$ does not
have any zero in the interval
$]T_1 , t_0]$. Without loss of generality, we can suppose that $t_0$ is the
first zero $> T_1$ of $v_1$.\\
From the indentities
$$\begin{array}{rcl}
[p_1(t)v_1'v_2]' & = & v_2[p_1(t)v_1']' + p_1(t)v_1'v_2' \\
\mbox{} & = & q_1(t)v_1v_2 - m_1(t)v_1v_2 + p_1(t)v_1'v_2'
\end{array}$$
and
$$\begin{array}{rcl}
[p_2(t)v_1v_2']' & = & v_1[p_2(t)v_2']' + p_2(t)v_1'v_2' \\
\mbox{} & = & q_2(t)v_1v_2 - m_2(t)v_1v_2 + p_2(t)v_1'v_2'
\end{array}$$
it follows that
$$[p_1(t)v_1'v_2 - p_2(t)v_1v_2']' =
[q_1(t) - q_2(t)]v_1v_2 -
[m_1(t) - m_2(t)]v_1v_2 + [p_1(t) - p_2(t)]v_1'v_2' \; .$$
Hence 
$$
\begin{array}{rcl}
\big(\frac{v_1}{v_2}[p_1(t)v_1'v_2 - p_2(t)v_1v_2']\big)' & = 
& \frac{v_1}{v_2}([q_1(t) - q_2(t)]v_1v_2 \\
\mbox{} & - & [m_1(t) - m_2(t)]v_1v_2 \\
\mbox{} & + & [p_1(t) - p_2(t)]v_1'v_2') \\
\mbox{} & + & \frac{v_1'v_2 - v_1v_2'}{v_2^2}[p_1(t)v_1'v_2 -
p_2(t)v_1v_2'] \; .
\end{array}
$$
So we establish the following formula 
\begin{equation}\label{eq:2.2}
\begin{array}{rcl}
\big(\frac{v_1}{v_2}[p_1(t)v_1'v_2 - p_2(t)v_1v_2']\big)' & = &   
[q_1(t) - q_2(t)]v_1^2\\
\mbox{} & - & [m_1(t) - m_2(t)]v_1^2\\
\mbox{} & + & [p_1(t) - p_2(t)]{v_1'}^2\\
\mbox{} & + & p_2(t)[v_1' - v_2'\frac{v_1}{v_2}]^2
\end{array} 
\end{equation} 
everywhere in the interval $]T_1 , t_0]$.\\
Now, integrating (\ref{eq:2.2}) from $T_1$ to $t_0$, the
first member leads to an integral equal to zero and the second one
to an integral strictly positive: a contradiction.\\
Similarly, we prove that $v_2$ must have a zero in $]T_1 , t_0[$ when at
least one of the inequalities $p_1 \geq p_2$, $q_1 \geq q_2$, $m_1 \leq m_2$
is strict. Q. E. D.
\begin{proper}
\label{proper:3}
$(i)\, \psi_1'(a) < 0$ for every $a \in B_1 \cap \dR^+$,\\
$(ii)\, \psi_1'(a) > 0$ for every $a \in B_1 \cap \dR^-$.
\end{proper}

{\bf Proof.} 
The proof is very similar to that of the third property of $\varphi$ (see
\cite{Al1} or \cite{Al2}). We point out that, instead of Sturm's comparison
theorem, one should use lemma \ref{lem:2.4} above. Q. E. D.
\vspace{3mm}
\par
The following property will be very useful in the fourth section of
this work. It concerns the
behaviour of $\psi_1$ at $\pm \infty$.
\begin{proper}
\label{proper:4}
(i) $\lim\limits_{a\longrightarrow +\infty}\psi_1(a) =
\alpha_{T_1}^{>}$ ,\\
(ii) $\lim\limits_{a\longrightarrow -\infty}\psi_1(a) = \alpha_{T_1}^<$.
\end{proper}

To prove this property, we need the following two lemmas:
\begin{lem}
\label{lem:2.5}
(i) $\forall a \in \dR_+ , \psi_1(a) \geq \alpha_{T_1}^{>}$ and if,
moreover, $a \in B_1$, then we have the strict inequality.\\
(ii) $\forall a \in \dR_- , \psi_1(a) \geq \alpha_{T_1}^{<}$ and if,
moreover, $a \in B_1$, then we have the strict inequality.
\end{lem}

\begin{lem}\label{lem:2.6}
$\forall a \in \dR , \psi_1(a) \leq \varphi_{a}(T_1)$. 
Moreover, if $(a , T_1) \in A$, then one has the strict inequality.
\end{lem}

{\bf Proof of property \ref{proper:4}} 
(i) Let $a \geq 0$. From the last two lemmas, we deduce that 
$$\alpha_{T_1}^{>} \leq \psi_1(a) \leq \varphi_a(T_1) \; .$$
Passing to the limit as
$a \longrightarrow +\infty$ and using the fact that
$\lim\limits_{a \longrightarrow +\infty}\varphi_{a}(T_1) = \alpha_{T_1}^{>}$
(see property 4 of $\varphi$ above), we conclude the proof of $(i)$.\\
(ii) The proof of (ii) is similar to that of (i). Q. E. D.
\vspace{3mm}
\par
{\bf Proof of lemma \ref{lem:2.5}.} 
Let us prove the first assertion $(i)$ (the second one is proved similarly).
Let $a \geq 0$ and suppose by contradiction that $\psi_1(a) < \alpha_{T_1}^>$.
Then $m$ is $\leq 0$ on $[T_1 , \psi_1(a)]$. Comparing on this interval our
equation (\ref{eq:1.6}) with the equation $-p_1 v'' = 0.v$,
it follows from lemma
\ref{lem:2.4} that any solution $v$ of the latter must have at least
one zero in $]T_1 , \psi_1(a)]$, which is clearly false.\\
By the same arguments, we prove that if moreover $a \in B_1$, then the
strict inequality $\psi_1(a) > \alpha_{T_1}^{>}$ holds true. Q. E. D.
\vspace{3mm}
\par
{\bf Proof of lemma \ref{lem:2.6}.} 
Let $a \in \dR$ be such that $(a , T_1) \in A$.
Applying Sturm's separation theorem as it is given for instance in
ch. 11 of \cite{Ha} to the functions $u(.;a,T_1)$ and
$v(.;a,T_1)$, it follows that $v(.;a,T_1)$ has at least one zero in the
interval $]T_1 , \varphi_{a}(T_1)[$ and,
hence, $\psi_1(a) < \varphi_{a}(T_1)$.\\
If $(a , T_1) \not\in A$, then it is clear by the definition of $A$ below that
$\varphi_{a}(T_1) = +\infty$ and, hence, $\psi_1(a) \leq \varphi_{a}(T_1)$.
Lemma \ref{lem:2.6} is proved. Q. E. D.
\begin{proper}
\label{proper:5}
$\lim\limits_{a \longrightarrow 0}\psi_{1}(a) = +\infty$.
\end{proper}

{\bf Proof.} 
Let us treat the case where $a \geq 0$ (the case
$a \leq 0$ can be reduced to this one by considering $-m$).\\
Let $R > 0$. It is clear that there exists $a_{R} > 0$ such that
$am(t) \leq p_{1}(\pi/R)^{2}$ for every
$0 \leq a \leq a_{R}$, where $p_{1}$ is defined at the beginning of this
section.\\
Let us compare, on $[T_1,T_1+R]$, the equations (\ref{eq:1.6}) and
\begin{equation}
-p_{1}v'' + 0.v = p_{1}(\pi/R)^{2}v.
\label{eq:2.3}
\end{equation}
The function $v(t) := \sin(\pi\frac {t-T_1}{R} + \frac {\pi}{2})$
is clearly a solution of (\ref{eq:2.3}) satisfying $v'(T_1) = 0$
and whose first zero $> T_1$ is exactly the point
$T_1 + \frac {R}{2}$. Using lemma \ref{lem:2.4}, this implies that 
$\psi_1(a) \geq T_1 + \frac {R}{2}$. Property \ref{proper:5} is then
proved since
$R$ is arbitrary. Q. E. D.
\begin{cons}
\label{cons:2.7}
From the remark \ref{rem:2.3} and the properties of $\psi_1$
presented above,
we deduce the following results concerning the zero-function $\psi_2$:
\vspace{2mm}\\
\begin{tabular}{cl}
1) & $B_2$ is an open subset of $\dR$ and
$\psi_2$ : $B_2 \longmapsto \dR$ is a $C^1$-function,\\
2) & $\psi_2$ is $\nearrow$ (resp. $\searrow$) in $\dR_+$ (resp.
$\dR_-$),
strictly in $B_2 \cap \dR_+$ (resp. $B_2 \cap \dR_-$),\\
3) & $\psi_{2}'(a)$ is $> 0$ if $a \in B_2 \cap \dR_+$ and
$< 0$ if $a \in B_2 \cap \dR_-$,\\
4) & $\lim\limits_{a \longrightarrow +\infty}\psi_{2}(a) =
\sup \{ t<T_2 : m(t)>0 \}$ and , similarly,\\
 & $\lim\limits_{a \longrightarrow -\infty}\psi_{2}(a) =
\sup \{ t<T_2 : m(t)<0 \}$,\\
5) & $\lim\limits_{a \longrightarrow 0}\psi_{2}(a) = -\infty$.
\end{tabular}
\end{cons}

Before closing this section, note that the restrictions of
\vspace{2mm}\\
\begin{tabular}{ccl}
$\varphi_a$ & to & $\{ s \in [T_1 , T_2] : \varphi_a(s) \in [T_1 , T_2]
\}$,\\
$\psi_1$ & to & $\{ a \in \dR : \psi_1(a) \in [T_1 , T_2] \}$ and \\
$\psi_2$ & to & $\{ a \in \dR : \psi_2(a) \in [T_1 , T_2] \}$,
\end{tabular}
\vspace{2mm}\\
clearly, do not depend on the extensions of the coefficients of $L$ and the
weight-function $m(t)$ introduced at the beginning of this section.\\
On the other hand, since $\varphi$, $\psi_1$ and $\psi_2$ depend on
the weight-function $m(t)$, we will note them below by
$\varphi^m$, $\psi_{1}^{m}$ and $\psi_{2}^{m}$ respectively.\\
Finally we denote by $\lambda_{1}^{m}$ (resp.
$\lambda_{-1}^{m}$) the principal positive (resp. negative) eigenvalue
of the problem (\ref{eq:1.4}) when $m^+(t)$ (resp. $m^-(t)$) $\not\equiv 0$
in $[T_1 , T_2]$.
\begin{rem}
\label{rem:2.8}
{\rm When $q \equiv 0$ and $m$ changes sign in $[T_1 , T_2]$,
it is well known that (at least) one of these
principal eigenvalues is equal to zero and the ``sign'' of the other one
depends on the ``average of the weight $m$'' on the interval $[T_1 , T_2]$.
More precisely, in this case,
\vspace{2mm}\\
\begin{tabular}{ccl}
{\bf .)} & if & $\int_{T_1}^{T_2} m(t)dt = 0$, then $\lambda_{-1}^{m} = 0
= \lambda_{1}^{m}$,\\
{\bf .)} & if & $\int_{T_1}^{T_2} m(t)dt > 0$, then $\lambda_{-1}^{m} < 0$
and $\lambda_{1}^{m} = 0$,\\
{\bf .)} & if & $\int_{T_1}^{T_2} m(t)dt < 0$, then $\lambda_{-1}^{m} = 0$
and $\lambda_{1}^{m} > 0$
\end{tabular}
\vspace{2mm}\\
(see, e.g., \cite{B-L} or \cite{S}).\par
When $q \geqs{\not \equiv} 0$ and $m$ changes sign in $[T_1 , T_2]$,
it is clear that $0$ can not be an
eigenvalue of the problem (\ref{eq:1.4}) and, moreover,
$\lambda_{-1}^m < 0$ and $\lambda_1^m > 0$ (see, e.g. \cite{H} for more
details).
}
\end{rem}

\section{Nonlinear problem with weights}

This section is devoted to the study of the description of the Fu\v cik
spectrum $\Sigma$ in the Neumann case. In the first paragraph, we
consider the problem with two weight-functions (\ref{eq:1.1}) and in the
second one, we investigate the particular case where one weight is
considered, i.e. the problem (\ref{eq:1.3}).

\subsection{Two-weights problem}
\setcounter{equation}{0}

Consider in this paragraph the semilinear problem (\ref{eq:1.1}), 
where $L$, $m$ and $n$ are defined as before. 
From standard ODE results it follows that the zeros of any nontrivial
solution of (\ref{eq:1.1}) are isolated.
Hence, we can classify these solutions according to their numbers of zeros in
the interval $]T_1 , T_2[$. This leads us to the following description
of the spectrum $\Sigma$:
$$\Sigma = \bigcup_{k = 0}^{+\infty} (C_k^> \cup C_k^<) \; ,$$
where $C_k^>$ (resp. $C_k^<$) is the set of those $(a,b)$ in $\dR^2$
such that (\ref{eq:1.1}) has a nontrivial solution $u$ with
exactly $k$-zeros in
$]T_1 , T_2[$ and ending positively (resp. negatively),
i.e. such that: $u(T_2) > 0$ (resp. $u(T_2) < 0$).\\
On the other hand, if $(a,b) \in C_k^>$ (resp. $C_k^<$), then all the
solutions of (\ref{eq:1.1}) ending positively (resp. negatively)
are multiple one of the other. This allows us to describe
$C_k^>$ and $C_k^<$, $k \geq 1$, as follows:
$$\begin{array}{ccl}
C_1^> & := & \{ (a,b) : \psi_1^n(b) = \psi_2^m(a) \} \, ,\\
C_2^> & := & \{ (a,b) : \varphi_b^n[\psi_1^m(a)] = \psi_2^m(a) \} \, ,\\
C_3^> & := & \{ (a,b) : \varphi_b^n(\varphi_a^m[\psi_1^n(b)]) = \psi_2^m(a) \} \, ,\\ 
C_4^> & := & \{ (a,b) : \varphi_b^n(\varphi_a^m[\varphi_b^n(\psi_1^m(a))]) = \psi_2^m(a) \} \; \ldots
\end{array}$$
and 
$$\begin{array}{ccl}
C_1^< & := & \{ (a,b) : \psi_1^m(a) = \psi_2^n(b) \} \, ,\\
C_2^< & := & \{ (a,b) : \varphi_a^m[\psi_1^n(b)] = \psi_2^n(b) \} \, ,\\
C_3^< & := & \{ (a,b) : \varphi_a^m(\varphi_b^n[\psi_1^m(a)]) = \psi_2^n(b) \} \, ,\\ 
C_4^< & := & \{ (a,b) : \varphi_a^m(\varphi_b^n[\varphi_a^m(\psi_1^n(b))]) = \psi_2^n(b) \} \; \ldots
\end{array}$$
Concerning the sets $C_0^>$ and $C_0^<$, it is clear that
$C_0^>$ (resp. $C_0^<$) is made of one or two of the lines
$\{\lambda_1^m\}\times\dR$ and $\{\lambda_{-1}^m\}\times\dR$ (resp.
$\dR\times\{\lambda_1^n\}$ and $\dR\times\{\lambda_{-1}^n\}$): it depends on
the sign of $m(t)$ (resp. $n(t)$). More precisely,\\
{\bf ./} If $m \geqs{\not\equiv} 0$ everywhere in
$[T_1 , T_2]$, then $C_0^> = \{\lambda_1^m\} \times \dR$.\\
{\bf ./} If $m \leqs{\not\equiv} 0$ everywhere in $[T_1 , T_2]$, then
$C_0^> = \{\lambda_{-1}^m\} \times \dR$.
\vspace{2mm}\\
{\bf ./} If $m$ changes sign in $[T_1 , T_2]$, i.e. if
$m^+(t) \not\equiv 0$ and $m^-(t) \not\equiv 0$ on $[T_1 , T_2]$, then
$C_0^> = (\{\lambda_{-1}^m\} \times \dR) \cup (\{\lambda_1^m\} \times
\dR)$.
\vspace{2mm}\\
Similarly,
\vspace{2mm}\\
.) $C_0^< = \dR \times \{\lambda_1^n\}$ if $n \geqs{\not\equiv} 0$
everywhere in
$[T_1 , T_2]$,\\
.) $C_0^< = \dR \times \{\lambda_{-1}^n\}$ if $n \leqs{\not\equiv} 0$
everywhere in
$[T_1 , T_2]$,\\
.) $C_0^< = (\dR \times \{\lambda_{-1}^n\}) \cup 
(\dR \times \{\lambda_1^n\})$ if $n$ changes sign in $[T_1 , T_2]$.
\vspace{3mm}
\par
From now on, we will denote by $\Sigma^*$ the set $\Sigma$ without
these trivial lines, i.e. 
$$\Sigma^* := \bigcup_{k = 1}^{+\infty} (C_k^> \cup C_k^<) \; .$$
\begin{rem}
\label{rem:3.1}
{\rm If for instance $m^+(t) \not \equiv 0$ and $n^+(t) \not \equiv 0$
in $[T_1 , T_2]$, then it follows
from the monotonicity properties of the zero-functions that the
intersection of $\Sigma$ with
$\dR^+ \times \dR^+$ is contained in the quadrant
$]\lambda_1^m , +\infty[ \times ]\lambda_1^n , +\infty[$. 
In fact, we can say more, namely this
inclusion is {\bf always strict}, i.e. there exists
$\epsilon > 0$ such that 
$\Sigma^* \cap (\dR^+ \times \dR^+) \subset ]\lambda_1^m + \epsilon ,
+\infty[ \times
]\lambda_1^n + \epsilon , +\infty[$ (see next section for more
details):
This is the most important difference between
the Dirichlet case and the Neumann case which is the subject of this work.
Indeed, it was proved in \cite{Al2} that this 
inclusion is strict in the Dirichlet case {\bf if and only if} 
$m^+(t)$ and $n^+(t)$ both have compact support in the 
interval $]T_1 , T_2[$.
\vspace{3mm}
\par
Now let us return to our Neumann case. If $m^+(t) \equiv 0$ and $n^+(t) \equiv 0$
in $[T_1 , T_2]$, then $\Sigma^*$ does not contain any point of
$\dR^+ \times \dR^+$.\\
We have similar results for the other three quadrants of
$\dR^2$: For instance, if $m(t)$ and $n(t)$ both change sign,
then $\Sigma^*$ is contained strictly in the four quadrants
$$(]\lambda_1^m , +\infty[ \times ]\lambda_1^n , +\infty[) \cup
(]-\infty , \lambda_{-1}^m[ \times ]-\infty , \lambda_{-1}^n[)$$ 
$$\cup
(]-\infty , \lambda_{-1}^m[ \times ]\lambda_1^n , +\infty[) \cup
(]\lambda_1^m , +\infty[ \times ]-\infty , \lambda_{-1}^n[).$$ 
When one of the weight-functions $m(t)$ or $n(t)$ changes sign, the other one
being with a constant sign in $]T_1 , T_2[$,
$\Sigma^*$ is contained strictly in two quadrants.\\
When these two weights do not change sign, 
$\Sigma^*$ is contained strictly in one quadrant.
}
\end{rem}

We will note below:
$$\begin{array}{ccl}
\Psi_1^>(a,b) & := & \psi_1^n(b) \, ,\\
\Psi_2^>(a,b) & := & \varphi_b^n[\psi_1^m(a)] \, ,\\
\Psi_3^>(a,b) & := & \varphi_b^n(\varphi_a^m[\psi_1^n(b)]) \, ,\\ 
\Psi_4^>(a,b) & := & \varphi_b^n(\varphi_a^m[\varphi_b^n(\psi_1^m(a))]) \,
\ldots
\end{array}$$
and, similarly, 
$$\begin{array}{ccl}
\Psi_1^<(a,b) & := & \psi_1^m(a) \, ,\\
\Psi_2^<(a,b) & := & \varphi_a^m[\psi_1^n(b)] \, ,\\
\Psi_3^<(a,b) & := & \varphi_a^m(\varphi_b^n[\psi_1^m(a)]) \, ,\\ 
\Psi_4^<(a,b) & := & \varphi_a^m(\varphi_b^n[\varphi_a^m(\psi_1^n(b))]) \,
\ldots
\end{array}$$
So we have:
$$\begin{array}{ccl}
C_{k}^{>} & = & \{(a,b) \in \dR^2 : \Psi_{k}^{>}(a,b) = \psi_2^m(a)\}\;
\mbox{and}\\
C_{k}^{<} & = & \{(a,b) \in \dR^2 : \Psi_{k}^{<}(a,b) = \psi_2^n(b)\}
\end{array}
$$
for every $k \geq 1$.
\vspace{3mm}
\par 
Taking into account of the fact that the zero-functions
$\varphi_.(T_1)$ and $\psi_1(.)$ have similar
properties, as it is shown in the preceding section,
and proceeding as in the Dirichlet case (\cite{Al1}, \cite{Al2} section 2),
we conclude that the properties of $\Sigma$ in our Neumann case are
rather comparable with those of $\Sigma$ in the Dirichlet case.\\
We point out that the asymptotic behaviour of $\Sigma$ represents an
exception of these properties (see remark \ref{rem:3.1}). This is the
principal reason for which we will devote the following section to this
fact.\\
So we present below the properties of the Fu\v cik spectrum associated to
(\ref{eq:1.1}). We refer to \cite{Al1} for detailed proofs and \cite{Al2}
for the proofs of analogous results in the Dirichlet case.
\vspace{3mm}
\par 
The first theorem shows that $\Sigma^*$ is made,
in each quadrant of $\dR^2$, of hyperbolic like curves of class $C^1$.
\begin{theo}
\label{theo:1} 
Let $k \geq 1$ be such that $C_{k}^{>} \cap (\dR_+ \times \dR_+)$
is nonempty. Then there exist $\alpha_{k}^{>} \geq \lambda_{1}^{m}$,
$\beta_{k}^{>} \geq \lambda_{1}^{n}$ and a $C^{1}$-diffeomorphism
$f_{k}^{>}$ strictly $\searrow$ defined from
$]\alpha_{k}^{>}, + \infty[$ to $]\beta_{k}^{>}, +\infty[$ such that 
   $$
   C_{k}^{>} \cap (\dR_+ \times \dR_+) =
   \{(a,f_{k}^{>}(a)): a \in ]\alpha_{k}^{>}, + \infty[\}.
   $$
Similarly for $C_{k}^{<}$.\\
We have similar results for the other three quadrants of
$\dR^2$. If for instance $C_{k}^{<} \cap (\dR_- \times \dR_+) \not
\equiv \emptyset$, then there exist $\gamma_k^< \leq \lambda_{-1}^m$,
$\eta_k^< \geq \lambda_1^n$ and a $C^{1}$-diffeomorphism $g_k^<$
strictly $\nearrow$ defined from $]-\infty , \gamma_k^<[$ to
$]\eta_k^< , +\infty[$ such that
$C_{k}^{<} \cap (\dR_- \times \dR_+) = \{ (a,g_k^<(a)) : a \in ]-\infty , 
\gamma_k^<[ \}$. 
\end{theo}

In the following two theorems, we give some
informations on the ``exact number'' of these curves in each
quadrant of $\dR^2$ and their positions with respect to each other.
\begin{theo}
\label{theo:2} 
Let $k \geq 1$ and $\epsilon_1$, $\epsilon_2 \in \{ +,- \}$.
Then the following two assertions are equivalent:\\
(i) $C_{k}^{>} \cap (\dR_{\epsilon_1} \times \dR_{\epsilon_2})$ {\bf and}
$C_{k}^{<} \cap (\dR_{\epsilon_1} \times \dR_{\epsilon_2})$
are both nonempty,\\
(ii) $C_{k+1}^{>} \cap (\dR_{\epsilon_1} \times \dR_{\epsilon_2})$ ({\bf 
or} 
$C_{k+1}^{<} \cap (\dR_{\epsilon_1} \times \dR_{\epsilon_2})$)
is nonempty.\\
On the other hand, both $C_k^>$ and $C_k^<$ are strictly above
(resp. below)
$C_{k+1}^{>}$ (or $C_{k+1}^{<}$) in the quadrants
$\dR_+ \times \dR_+$ and $\dR_- \times \dR_+$ (resp. $\dR_+ \times \dR_-$
and $\dR_- \times \dR_-$).
\end{theo}

\begin{theo}
\label{theo:3}
Let $\epsilon_1 , \epsilon_2 \in \{ +,- \}$ and suppose that
$m^{\epsilon_1}(t) \not \equiv 0$ and $n^{\epsilon_2}(t) \not \equiv 0$.
Then one, at least, of the two intersections
$C_{1}^{>} \cap (\dR_{\epsilon_1} \times \dR_{\epsilon_2})$
or $C_{1}^{<} \cap (\dR_{\epsilon_1} \times \dR_{\epsilon_2})$ is nonempty.
\end{theo}

Combining theorems \ref{theo:2} and \ref{theo:3} yields the following
\begin{cons}\label{cons:3.2}  
Let $\epsilon_1 , \epsilon_2 \in \{ +,- \}$ and suppose
$m^{\epsilon_1}(t) \not \equiv 0$ and $n^{\epsilon_2}(t) \not \equiv 0$.
Then\par
either\\
(i) For every integer $k$, the intersections
$C_{k}^{>} \cap (\dR_{\epsilon_1} \times \dR_{\epsilon_2})$ and
$C_{k}^{<} \cap (\dR_{\epsilon_1} \times \dR_{\epsilon_2})$ are both
nonempty,\par
or\\
(ii) There exists $k_0$ for which
$C_{k}^{>} \cap (\dR_{\epsilon_1} \times \dR_{\epsilon_2})$ and
$C_{k}^{<} \cap (\dR_{\epsilon_1} \times \dR_{\epsilon_2})$ are
both nonempty for every $k \leq k_0$, one of the two intersections
$C_{k_0+1}^{>} \cap (\dR_{\epsilon_1} \times \dR_{\epsilon_2})$ or
$C_{k_0+1}^{<} \cap (\dR_{\epsilon_1} \times \dR_{\epsilon_2})$ is nonempty,
the other one being empty, and
$C_{k}^{>} \cap (\dR_{\epsilon_1} \times \dR_{\epsilon_2})$ and
$C_{k}^{<} \cap (\dR_{\epsilon_1} \times \dR_{\epsilon_2})$ are
both empty for every $k \geq k_0+2$.
\par
In particular, we conclude that if both $m(t)$ and
$n(t)$ change sign, then the spectrum $\Sigma^*$
consists, in each quadrant of $\dR^2$, of {\bf a (non-zero) odd
or infinite number} of hyperbolic like curves $C_k^> , C_k^< (k \geq 1)$.
\end{cons}

\begin{rem}
\label{rem:3.3}  
{\rm Note that certain curves may be double and are then counted for two.}
\end{rem}
 
The following theorem gives a sufficient and (almost)
necessary condition on $m(t)$ and $n(t)$ in order to have, in a given
quadrant of $\dR^2$, an infinite number of curves of $\Sigma^*$ 
(see \cite{Al1} for a more general condition which is necessary and
sufficient). 
\begin{theo}
\label{theo:4}
Let $\epsilon_1 , \epsilon_2 \in \{ +,- \}$ and suppose that
$m^{\epsilon_1}(t).n^{\epsilon_2}(t) \not \equiv 0$ in $[T_1 , T_2]$.
Then the two intersections
$C_{k}^{>} \cap (\dR_{\epsilon_1} \times \dR_{\epsilon_2})$ and
$C_{k}^{<} \cap (\dR_{\epsilon_1} \times \dR_{\epsilon_2})$ are
both nonempty for every $k \geq 1$.\\
Moreover, if at least one of the two graphes of $m^{\epsilon_1}(t)$
or $n^{\epsilon_2}(t)$ is made of a finite number of bumps,
then the converse is true.
\end{theo}

The last result shows that all situations may happen concerning
the number of curves of $\Sigma^*$ contained in each quadrant of $\dR^2$. 
\begin{theo}
\label{theo:5}
Let $p, q, r, s \in \dN \cup \{ +\infty \}$. Then there exist two
continuous weight-functions $m(t)$ and $n(t)$ on $[T_1 , T_2]$ such that
the intersection between the spectrum $\Sigma^*$ associated to the
problem (\ref{eq:1.1}) and $\dR_+ \times \dR_+$
(resp. $\dR_+ \times \dR_-$, $\dR_- \times \dR_-$, $\dR_- \times \dR_+$)
exactly contains $2p+1$ (resp. $2q+1$, $2r+1$, $2s+1$)
hyperbolic like curves.
\end{theo}

\subsection{Particular case where $m(t) \equiv n(t)$}

Now, let us consider the problem (\ref{eq:1.3}) where $L$ is defined as
before and $m \in C[T_1 , T_2]$,
$m(t) \not \equiv 0$. Making $a = b$ in (\ref{eq:1.3}),
we find the eigenvalue problem (\ref{eq:1.4}) and, hence,
we conclude that
$\Sigma$ intersected by the diagonal ``$y=x$'' consists of the sequence
of points $(\lambda_k^m , \lambda_k^m)$, $k \geq 1$,
if $m^+(t) \not \equiv 0$ and
$(\lambda_{-k}^m , \lambda_{-k}^m)$, $k \geq 1$, if $m^-(t) \not \equiv 0$.
Moreover, it is clear that
$(\lambda_k^m , \lambda_k^m) \in (C_{k-1}^> \cap C_{k-1}^<)
\cap \dR_+^2$ and
$(\lambda_{-k}^m , \lambda_{-k}^m) \in (C_{k-1}^> \cap C_{k-1}^<)
\cap \dR_-^2$ for every $k \geq 1$. This implies that:\\
(i) if $m^+(t) \not \equiv 0$, then all the intersections of
$C_k^>$ and of $C_k^<$ with $\dR_+^2$ are nonempty,\\
(ii) if $m^-(t) \not \equiv 0$, then all the intersections of
$C_k^>$ and of $C_k^<$ with $\dR_-^2$ are nonempty.\\
In particular, if $m(t)$ does not change sign in the interval
$[T_1 , T_2]$ and if $m \geqs{\not \equiv} 0$ (resp. $m \leqs{\not \equiv}
0$), then $\Sigma^*$ contains an infinite number of curves in the quadrant
$\dR_+^2$ (resp. $\dR_-^2$).\par
From now on, we suppose that the weight-function
{\bf $m(t)$ changes sign}, i.e. $m^+$ and $m^-$ are both
$\not \equiv 0$ in $[T_1 , T_2]$.\\
Using theorem \ref{theo:1}, it follows from the points (i) and (ii)
above that $\Sigma^*$ contains an infinite
number of curves in each one of the quadrants $\dR_+^2$ and $\dR_-^2$.
Moreover, using the fact that $\Sigma^*$ is symmetric with respect
to the diagonal ``$y=x$'',
it follows from the consequence \ref{cons:3.2} that:
either
(i) $\Sigma^*$ contains a (non-zero) odd number of hyperbolic like
curves in each one of the two quadrants $\dR_+ \times \dR_-$ and
$\dR_- \times \dR_+$, or (ii) there exists an infinite number of
hyperbolic like curves in each one of these two quadrants.
\vspace{3mm}
\par 
We close this section by the following theorem which shows that the
``exact number'' of these additional curves in $\dR_+ \times \dR_-$ and
in $\dR_- \times \dR_+$ depends on the ``number of changes of sign''
of $m(t)$ in the interval $]T_1 , T_2[$.\\ 
Let us first precise this notion of ``number of changes of sign''.
\begin{defi}
\label{defi:3.4} 
   {\rm
   {\bf (1)} Let $s \in ]T_{1},T_{2}[$. We say that $s$ is a {\sl
   simple}
   point of change of sign of $m$ when there exist
   $T_{1} < s' \leq s$ and
    $\epsilon_{0} > 0$ such that: either (i) $m \leqs{\not \equiv} 0$ 
    on $]s'-\epsilon,s'[ , \forall 0 < \epsilon < \epsilon_{0}$, $m 
    \equiv 0$ on $[s',s]$ and $m \geqs{\not \equiv} 0$ on
    $]s,s+\epsilon[ , \forall 0 < \epsilon < \epsilon_0$, or
    (ii) $m \geqs{\not \equiv} 0$ on
    $]s'-\epsilon,s'[$ for every $0 <  \epsilon < \epsilon_{0}$, 
    $m\equiv 0$ on $[s',s]$ and $m \leqs{\not \equiv} 0$ on 
    $]s,s+\epsilon[$ for every $0 < \epsilon < \epsilon_{0}$.\\
    {\bf (2)} Let $s \in [T_{1},T_{2}]$. We say that $s$ is a {\sl
    multiple}
    point of change of sign of $m$ if either (i) $s > T_1$ and $m^+$ and $m^-$
    are both $\not \equiv 0$ on $]s-\epsilon,s[\, \cap\, 
    [T_1 , T_2]$ for every $\epsilon > 0$, or (ii) $s < T_2$ and $m^+$ and
    $m^-$ are $\not \equiv 0$ on $]s,s+\epsilon[\, \cap\, [T_1 , T_2]$
    for any $\epsilon > 0$.\\
   {\bf (3)} If $m$ has only simple points of change of sign in
    $]T_1 , T_2[$ and if $N \in \dN$ is the exact
    number of these points, then we say that the number of changes of sign
    of $m$ is equal to $N$. If $m$ has, at least, one multiple point
    of change of sign in $[T_1 , T_2]$, then we say that the number of changes
    of sign of $m$ is $+\infty$.
}
\end{defi}

Hence we have the following
\begin{theo}
\label{theo:6}
Suppose that $m$ changes sign in the interval $]T_1 , T_2[$.
Let $N \in \{ 1, 2, ..., + \infty \}$ be the number of changes of sign of $m$.
Then the spectrum $\Sigma^*$ exactly contains $(2N-1)$-hyperbolic like
curves in each one of the two quadrants $\dR_+ \times \dR_-$ and
$\dR_- \times \dR_+$.
\end{theo}

\section{Asymptotic behaviour of the first curves of $\Sigma$}
\setcounter{equation}{0}

This section is devoted to the study of the asymptotic behaviour of the
first curves $C_1^> \cap (\dR_{\epsilon_1} \times \dR_{\epsilon_2})$ and
$C_1^< \cap (\dR_{\epsilon_1} \times \dR_{\epsilon_2})$, where
$\epsilon_1 , \epsilon_2 \in \{ +,- \}$, of the Fu\v cik spectrum
$\Sigma$ associated
to the problem (\ref{eq:1.3}). We will see that none of the first curves can
be asymptotic on any side to the trivial
lines of $\Sigma$.\\
We denote
\begin{eqnarray*}
    T_{1}^{>} & := & \inf \{ t \in ]T_{1}, T_{2}[ : m(t) > 0\},  \\
    T_{2}^{>} & := & \sup \{ t \in ]T_{1}, T_{2}[ : m(t) > 0\}
\end{eqnarray*}
when $m^+(t) \not \equiv 0$ in $[T_1 , T_2]$, and
\begin{eqnarray*}
    T_{1}^{<} & := & \inf \{ t \in ]T_{1}, T_{2}[ : m(t) < 0\},  \\
    T_{2}^{<} & := & \sup \{ t \in ]T_{1}, T_{2}[ : m(t) < 0\}
\end{eqnarray*}
when $m^-(t) \not \equiv 0$ in $[T_1 , T_2]$.\\
For every $\epsilon_1 , \epsilon_2 \in \{ +,- \}$, we set
\begin{eqnarray*}
    (C_1^>)_{(\epsilon_1 , \epsilon_2)} & := & 
    C_1^> \cap (\dR_{\epsilon_1} \times \dR_{\epsilon_2}),  \\
    (C_1^<)_{(\epsilon_1 , \epsilon_2)} & := & 
    C_1^< \cap (\dR_{\epsilon_1} \times \dR_{\epsilon_2}). 
\end{eqnarray*}
For example, $(C_1^<)_{(-,+)} := C_1^< \cap (\dR_- \times \dR_+)$.
\par
Finally, we denote
\begin{enumerate}
\item[(i)]
$\alpha_1^{>}$ (resp. $\alpha_1^{<}$) the real $> 0$ (resp.
$< 0$)
for which $\psi_2(\alpha_1^>) = T_1^>$ (resp. $\psi_2(\alpha_1^<) =
T_1^>$),

\item[(ii)]
$\beta_1^{>}$ (resp. $\beta_1^{<}$) the real $> 0$ (resp.
$< 0$)
for which $\psi_1(\beta_1^>) = T_2^>$ (resp. $\psi_1(\beta_1^<) =
T_2^>$)
\end{enumerate}
if $m^+(t) \not \equiv 0$ in $[T_1 , T_2]$, and
\begin{enumerate}
\item[(iii)]
$\alpha_2^{>}$ (resp. $\alpha_2^{<}$) the real $> 0$ (resp.
$< 0$)
for which $\psi_2(\alpha_2^>) = T_1^<$ (resp. $\psi_2(\alpha_2^<) =
T_1^<$),

\item[(iv)]
$\beta_2^{>}$ (resp. $\beta_2^{<}$) the real $> 0$ (resp. $<
0$)
for which $\psi_1(\beta_2^>) = T_2^<$ (resp. $\psi_1(\beta_2^<) =
T_2^<$)
\end{enumerate}
if $m^-(t) \not \equiv 0$ in $[T_1 , T_2]$.\\
Note that the unicity of these reals clearly follows from the strict
monotonicity of the zero-functions $\psi_1$ and $\psi_2$. 
\begin{theo}
\label{theo:7}
    \begin{enumerate}
        \item[(i)] If $(C_1^>)_{(+,+)} \neq \emptyset$, then 
        $(C_1^>)_{(+,+)}$ is asymptotic to the lines
        $\{\alpha_1^>\} \times \dR$ and $\dR \times \{\beta_1^>\}$.

        \item[(ii)] If $(C_1^>)_{(-,-)} \neq \emptyset$, then 
        $(C_1^>)_{(-,-)}$ is asymptotic to the lines
        $\{\alpha_2^<\} \times \dR$ and $\dR \times \{\beta_2^<\}$.
    
        \item[(iii)] If $(C_1^>)_{(+,-)} \neq \emptyset$, then 
        $(C_1^>)_{(+,-)}$ is asymptotic to the lines
        $\{\alpha_2^>\} \times \dR$ and $\dR \times \{\beta_1^<\}$.
    
        \item[(iv)] If $(C_1^>)_{(-,+)} \neq \emptyset$, then 
        $(C_1^>)_{(-,+)}$ is asymptotic to the lines
        $\{\beta_2^>\} \times \dR$ and $\dR \times \{\alpha_1^<\}$.
\vspace{1mm}\\
And, by symmetry of $C_1^>$ and $C_1^<$ with respect to the diagonal
``$y=x$'', we have also:
        \item[(i)'] If $(C_1^<)_{(+,+)} \neq \emptyset$, then 
        $(C_1^<)_{(+,+)}$ is asymptotic to the lines
        $\{\beta_1^>\} \times \dR$ and $\dR \times \{\alpha_1^>\}$.

        \item[(ii)'] If $(C_1^<)_{(-,-)} \neq \emptyset$, then 
        $(C_1^<)_{(-,-)}$ is asymptotic to the lines
        $\{\beta_2^<\} \times \dR$ and $\dR \times \{\alpha_2^<\}$.
    
        \item[(iii)'] If $(C_1^<)_{(+,-)} \neq \emptyset$, then 
        $(C_1^<)_{(+,-)}$ is asymptotic to the lines
        $\{\alpha_1^<\} \times \dR$ and $\dR \times \{\beta_2^>\}$.
    
        \item[(iv)'] If $(C_1^<)_{(-,+)} \neq \emptyset$, then
        $(C_1^<)_{(-,+)}$ is asymptotic to the lines
        $\{\beta_1^<\} \times \dR$ and $\dR \times \{\alpha_2^>\}$.
    \end{enumerate}
\end{theo}

{\bf Proof.} Let $(a,b) \in C_1$. Then $\psi_1(b) = \psi_2(a)$ if
$(a,b) \in C_1^>$, or $\psi_1(a) = \psi_2(b)$ if $(a,b) \in C_1^<$.\\
Passing to the limits as $a \longrightarrow \pm \infty$ 
and as $b \longrightarrow \pm \infty$ and using 
the property \ref{proper:4} above and the fact that
$\psi_1$ and $\psi_2$ are strictly monotonous, one can easily prove
all the assertions of theorem \ref{theo:7}. Q. E. D.
\begin{coro}
\label{coro:4.1}
    If $m(t)$ changes sign in $[T_1 , T_2]$, then there exists an
    $\epsilon > 0$ such that
    $$\Sigma^* \subset
    (]\lambda_1^m + \epsilon , +\infty[ \times
    ]\lambda_1^m + \epsilon , +\infty[) \cup
    (]-\infty , \lambda_{-1}^m - \epsilon[ \times
    ]-\infty , \lambda_{-1}^m - \epsilon[)$$
    $$\cup
    (]-\infty , \lambda_{-1}^m - \epsilon[ \times
    ]\lambda_1^m + \epsilon , +\infty[) \cup
    (]\lambda_1^m + \epsilon , +\infty[ \times
    ]-\infty , \lambda_{-1}^m - \epsilon[).$$
    If $m(t)$ keeps a constant sign in $[T_1 , T_2]$ (for instance $m
    \geqs{\not \equiv} 0$), then 
    $$\Sigma^* \subset\, 
    ]\lambda_1^m + \epsilon , +\infty[ \times
    ]\lambda_1^m + \epsilon , +\infty[$$
    for a certain $\epsilon > 0$.
\end{coro}

{\bf Proof.} Suppose that $m(t)$ changes sign in $[T_1 , T_2]$. We denote
\vspace{2mm}\\
\begin{tabular}{cr}
$\alpha^>$ & the real $> 0$ for which $\psi_2(\alpha^>) = T_1$,\\
$\alpha^<$ & the real $< 0$ for which $\psi_2(\alpha^<) = T_1$,\\
$\beta^>$ & the real $> 0$ for which $\psi_1(\beta^>) = T_2$,\\
$\beta^<$ & the real $< 0$ for which $\psi_1(\beta^<) = T_2$.
\end{tabular}
\vspace{2mm}\\
Hence, by the monotonicity of the zero-functions $\psi_1$ and $\psi_2$,
it is clear that:
\vspace{2mm}\\
$\alpha_1^> , \alpha_2^> \geq \alpha^> > \lambda_1^m$ ;\\
$\alpha_1^< , \alpha_2^< \leq \alpha^< < \lambda_{-1}^m$ ;\\
$\beta_1^> , \beta_2^> \geq \beta^> > \lambda_1^m$ and\\
$\beta_1^< , \beta_2^< \leq \beta^< < \lambda_{-1}^m$.
\vspace{2mm}\\
So, in this case, one has just to take
$\epsilon := \min \{ \alpha^> - \lambda_1^m , \lambda_{-1}^m - \alpha^< ,
\beta^> - \lambda_1^m , \lambda_{-1}^m - \beta^< \}$.\\
If $m(t)$ is $\geqs{\not \equiv} 0$, then we take,
for example, 
$\epsilon := \min \{ \alpha^> - \lambda_1^m , \beta^> - \lambda_1^m \}$.
Q. E. D.
\begin{rem}
\label{rem:4.2}
{\rm
    By the same arguments, one could prove that the spectrum
    $\Sigma^*$ associated to the two-weights problem (\ref{eq:1.1})
    is always at a strictly positive distance from the trivial lines
    of $\Sigma$. 
    In the Dirichlet case, it was proved in \cite{Al2} that
    this is true {\bf if and only if} each one of the two weight-functions
    $m(t)$ and $n(t)$ has a compact support in the interval
    $]T_1 , T_2[$.
}
\end{rem}


\begin{thebibliography}{xx}

\bibitem{Al1} M. ALIF, Spectre de Fu\v cik : probl\`eme avec poids en 
dimension un et quelques remarques
en dimension sup\'erieure, Ph. D. Thesis,
Universit\'e Libre de Bruxelles, 1999.

\bibitem{Al2} M. ALIF \& J.-P. GOSSEZ, On the Fu\v cik spectrum with
indefinite weights, Diff. Int. Equat., to appear.

\bibitem{B-L} K. J. BROWN \& S. S. LIN, On the existence of positive
eigenfunctions for an eigenvalue problem with an indefinite weight
function, Mat. Anal. and Appl. 75, 112-120 (1980).

\bibitem{Ca} J. CAMPOS, Espectro de Fu\v cik para operadores elipticos,
Ph. D. Thesis, Universidad de Granada, 1996.

\bibitem{Co} E. A. CODDINGTON \& N. LEVINSON, Theory of ordinary
differential equations, McGraw-Hill, New York, 1955.

\bibitem{Da} N. DANCER, On the Dirichlet problem for weakly nonlinear
elliptic partial differential equations, Proc.
Royal Soc., Edimb., 76(1977), 283-300.

\bibitem{Dr}  P. DRABEK, Solvability and bifurcations of nonlinear 
equations, Pitman Research Notes in Mathematics, 264 (1992).

\bibitem{Fu} S. FU\v CIK, Solvability of nonlinear equations and boundary
value problems, Reidel, Dordrecht, 1980.

\bibitem{Ha} P. HARTMAN, Ordinary differential equations, Wiley, 1964.

\bibitem{H} P. HESS, On the spectrum of elliptic operators with respect
to indefinite weights, Linear Algebra and its Applications, 84(1986),
99-109.

\bibitem{Ry} B. RYNNE, The Fu\v cik spectrum of general Sturm-Liouville
problems, J. Diff. Equat., to appear.

\bibitem{S} S. SENN, On a nonlinear elliptic eigenvalue problem with
Neumann boundary conditions, with an application to population genetics,
Comm. in Part. Diff. Eq., 8(11), 1199-1228 (1983).
\end{thebibliography}
\end{document}